# STRONG INVARIANCE PRINCIPLES FOR SEQUENTIAL BAHADUR–KIEFER AND VERVAAT ERROR PROCESSES OF LONG-RANGE DEPENDENT SEQUENCES

By Miklós Csörgő,[1] Barbara Szyszkowicz[1] and Lihong Wang

*Carleton University, Carleton University and Nanjing University*

In this paper we study strong approximations (invariance principles) of the sequential uniform and general Bahadur–Kiefer processes of long-range dependent sequences. We also investigate the strong and weak asymptotic behavior of the sequential Vervaat process, that is, the integrated sequential Bahadur–Kiefer process, properly normalized, as well as that of its deviation from its limiting process, the so-called Vervaat error process. It is well known that the Bahadur–Kiefer and the Vervaat error processes cannot converge weakly in the i.i.d. case. In contrast to this, we conclude that the Bahadur–Kiefer and Vervaat error processes, as well as their sequential versions, do converge weakly to a Dehling–Taqqu type limit process for certain long-range dependent sequences.

**1. Introduction.** Assume that we have a stationary long-range dependent sequence of standard Gaussian random variables, $\eta_1, \eta_2, \ldots, \eta_n, \ldots$, that is, the Gaussian sequence $\{\eta_n, n \geq 1\}$ with $E\eta_1 = 0$ and $E\eta_1^2 = 1$ is assumed to have a positive covariance function of the form

$$(1.1) \qquad \gamma(k) := E(\eta_1 \eta_{k+1}) = k^{-D} L(k), \qquad 0 < D < 1,$$

for large $k$, where $L(\cdot)$ is a slowly varying function at infinity in the sense that

$$\lim_{s \to \infty} \frac{L(st)}{L(s)} = 1 \qquad \text{for every } t \in (0, \infty).$$

Received November 2003; revised June 2005.
[1]Supported in part by NSERC Canada grants of Miklós Csörgő and Barbara Szyszkowicz at Carleton University, Ottawa.

Let $G$ be an arbitrary real-valued Borel measurable function on the real line $\mathbb{R}$, and consider the subordinate process

$$(1.2) \qquad X_n = G(\eta_n), \qquad n \geq 1,$$

with marginal distribution function $F(x) = P(X \leq x)$, $x \in \mathbb{R}$, where $X = G(\eta)$ and $\eta$ is a standard normal random variable.

Assumption (1.2) allows one to use the theory of nonlinear functionals of Gaussian processes. As in [20], we expand the function $I(X_n \leq x) - F(x) = I(G(\cdot) \leq x) - F(x)$ in Hermite polynomials, for any fixed $x \in \mathbb{R}$,

$$I(X_n \leq x) - F(x) = \sum_{l=\tau_x}^{\infty} c_l(x) H_l(\eta_n)/l!,$$

where

$$H_l(x) = (-1)^l e^{x^2/2} \frac{d^l}{dx^l} e^{-x^2/2}, \qquad l = 1, 2, \ldots, x \in \mathbb{R},$$

is the $l$th Hermite polynomial,

$$c_l(x) = \mathrm{E}\{[I(G(\eta) \leq x) - F(x)] H_l(\eta)\},$$

and $\tau_x$ for any $x \in \mathbb{R}$ is the index of the first nonzero coefficient in the expansion, called the Hermite rank of the function $I(G(\cdot) \leq x) - F(x)$. We note in passing that, by the Cauchy–Schwarz inequality,

$$|c_l(x)| \leq \mathrm{E}(H_l^2(\eta)) = l!.$$

As in [20], the Hermite rank of the class of functions $\{I(X_n \leq x) - F(x), x \in \mathbb{R}\}$ is defined by

$$(1.3) \qquad \tau = \min\{\tau_x : c_{\tau_x}(x) \neq 0 \text{ for some } x \in \mathbb{R}\},$$

that is, $\tau = \inf_x \tau_x$. If we assume that $F$ is continuous, then the induced sequence of random variables

$$(1.4) \qquad U_n = F(X_n) = F(G(\eta_n)), \qquad n \geq 1,$$

is a Uniform-$[0,1]$ random sequence. Consequently, for any fixed $y \in (0,1)$, the function $(I(U_n \leq y) - y) = (I(F(G(\cdot)) \leq y) - y)$ has the Hermite expansion

$$I(U_n \leq y) - y = \sum_{l=\tau}^{\infty} J_l(y) H_l(\eta_n)/l!,$$

where

$$J_l(y) = \mathrm{E}\{[I(F(G(\eta)) \leq y) - y] H_l(\eta)\}.$$

Obviously, $J_l(y) = c_l(Q(y))$ for any $y \in (0,1)$, where $Q$ is the quantile function of $F$, that is,

$$Q(y) = F^{-1}(y) = \inf\{x : F(x) = y\}, \qquad 0 < y \leq 1, Q(0) = Q(0+),$$

and, hence, the Hermite rank of the class of functions $\{I(U_n \leq y) - y, y \in (0,1)\}$ is also $\tau$.

Given chronologically ordered samples $X_1, \ldots, X_n$ and $U_1, \ldots, U_n$, $n \geq 1$, as in (1.2) and (1.4), respectively, their corresponding sequential empirical distribution functions are

$$\widehat{E_{[nt]}}(y) = \begin{cases} 0, & 0 \leq t < 1/n, \\ \dfrac{1}{[nt]} \sum_{i=1}^{[nt]} I(U_i \leq y), & 0 \leq y \leq 1, 1/n \leq t \leq 1, \end{cases}$$

and

$$\widehat{F_{[nt]}}(x) = \begin{cases} 0, & 0 \leq t < 1/n, \\ \dfrac{1}{[nt]} \sum_{i=1}^{[nt]} I(X_i \leq x), & -\infty < x < \infty, 1/n \leq t \leq 1. \end{cases}$$

Based on these functions, we define the sequential empirical quantile functions

$$\widehat{U_{[nt]}}(y) = \widehat{E_{[nt]}}^{-1}(y) = \inf\{s : \widehat{E_{[nt]}}(s) \geq y\}, \qquad 0 < y \leq 1,$$
$$\widehat{U_{[nt]}}(0) = \widehat{U_{[nt]}}(0+), \qquad 0 \leq t \leq 1,$$

and

$$\widehat{Q_{[nt]}}(y) = \widehat{F_{[nt]}}^{-1}(y) = \inf\{x : \widehat{F_{[nt]}}(x) \geq y\}, \qquad 0 < y \leq 1,$$
$$\widehat{Q_{[nt]}}(0) = \widehat{Q_{[nt]}}(0+), \qquad 0 \leq t \leq 1.$$

Now the corresponding sequential uniform and general empirical and quantile processes are defined by

$$\alpha_n(y,t) = d_n^{-1}[nt](\widehat{E_{[nt]}}(y) - y), \qquad 0 \leq y \leq 1, 0 \leq t \leq 1,$$
$$u_n(y,t) = d_n^{-1}[nt](y - \widehat{U_{[nt]}}(y)), \qquad 0 \leq y \leq 1, 0 \leq t \leq 1,$$
$$\beta_n(x,t) = d_n^{-1}[nt](\widehat{F_{[nt]}}(x) - F(x)), \qquad -\infty < x < \infty, 0 \leq t \leq 1,$$
$$\gamma_n(y,t) = d_n^{-1}[nt](Q(y) - \widehat{Q_{[nt]}}(y)), \qquad 0 < y < 1, 0 \leq t \leq 1,$$

where

(1.5) $$d_n^2 = n^{2-\tau D} L^\tau(n),$$

with $D$ of (1.1) so that $0 < D < 1/\tau$, where $\tau$ is defined in (1.3).

By Theorem 3.1 of [31], one arrives at

$$\text{Var}(n\widehat{F}_n(x)) \sim n^{2-\tau D} L^\tau(n) \frac{2c_\tau^2(x)}{\tau!(2-\tau D)(1-\tau D)} = O(d_n^2)$$

for each fixed $x \in \mathbb{R}$ as $n \to \infty$, where the symbol $\sim$ means asymptotic proportional equivalence. This explains the choice of $d_n$ as defined in (1.5) for defining the above sequential empirical and quantile processes.

Dehling and Taqqu [19, 20] studied the asymptotic properties of the sequential general empirical process $\beta_n(x,t)$. The following important two-parameter weak convergence theorem for $\beta_n(x,t)$ is due to Dehling and Taqqu [20] whose Theorem 1.1 reads as follows.

THEOREM A. *Let the stationary subordinate process $\{X_n, n \geq 1\}$ be as in (1.2) with $\tau$ as in (1.3), and let $d_n$ be as in (1.5). Then, as $n \to \infty$,*

$$\{\beta_n(x,t); -\infty \leq x \leq +\infty, 0 \leq t \leq 1\}$$

*converges weakly in $D[-\infty, +\infty] \times [0,1]$,*

*equipped with the sup-norm, to*

$$\left\{c_\tau(x)\sqrt{\frac{2}{(2-\tau D)(1-\tau D)}} Y_\tau(t); -\infty \leq x \leq +\infty, 0 \leq t \leq 1\right\},$$

$$0 < D < 1/\tau,$$

*where $Y_\tau(t)$ is $1/\tau!$ times a Hermite process of rank $\tau$, given for each $t \in [0,1]$ as a multiple Wiener–Itô–Dobrushin integral that is defined in* (1.7) *of* [20].

Thus, in Theorem A, $\tau! Y_\tau(t) =: Z_\tau(t)$ is a Hermite process of rank $\tau$, a self-similar, stationary increment process with self-similarity index $H = 1 - \tau D/2$, $0 < D < 1/\tau$, which, as shown in [21, 31, 33], can be represented as the multiple Wiener–Itô stochastic integral

$$Z_\tau(t) := \int_{\mathbb{R}^\tau} k_t(x_1, \ldots, x_\tau) \, dW(x_1) \cdots dW(x_\tau)$$

with respect to a standard Wiener process (Brownian motion) $W(x)$, $E(dW(x))^2 = dx$, where

$$k_t(x_1, \ldots, x_\tau) = K(\tau, D) \int_0^t \prod_{i=1}^\tau ((s-x_i)^+)^{-(D+1)/2} \, ds$$

and

$$K(\tau, D) = ((2-\tau D)(1-\tau D)/2) \left\{\frac{1}{\tau!}\left(\frac{1}{B(D, (1-D)/2)}\right)^\tau\right\}^{1/2},$$

with $B(\cdot,\cdot)$ denoting the beta-function.

For a general theory of multiple stochastic integration and Hermite processes, we refer to [27].

We also note that Dehling and Taqqu [19] obtained a functional law of the iterated logarithm as well for $\beta_n(x,t)$ in $D[-\infty,+\infty] \times [0,1]$.

REMARK 1.1. We recall (cf. [29]) that, in the i.i.d. case, the weak limit of $\beta_n(x,t)$ is a two-time parameter Gaussian process in $x$ and $t$, the so-called Kiefer process on account of the landmark Kiefer [25] paper, which is a Brownian bridge in $x$ and a Wiener process (Brownian motion) in $t$. The Dehling–Taqqu [20] limit in Theorem A differs greatly from the Kiefer process. Namely, it separates the variables in $x$ and $t$ in terms of being the product of a *deterministic function* in $x$ and a *stochastic process* in $t$ which is non-Gaussian when $\tau \geq 2$. In particular, if in (1.2) $G(x) = x$, then $\tau$ of (1.3) is equal to 1, and $Y_1(t) = Z_1(t)$ of Theorem A is a fractional Brownian motion with self-similarity index $H = 1 - D/2$, $0 < D < 1$. If $G(x) = x^2$, then $\tau = 2$, and $Y_2(t)$ of Theorem A equals $2^{-1}Z_2(t)$, where $Z_2(t)$ is non-Gaussian, has stationary increments and the same covariance as $Z_1(t)$ but with $H = 1 - D$, $0 < D < 1/2$. It is called the Rosenblatt process (cf. [31]). For details on the latter two examples, and on that of Hermite rank $\tau > 2$, we refer to pages 1770–1771 of [20].

Assuming that $F$ has a Lebesgue density function $f$ on $\mathbb{R}$, Csörgő and Mielniczuk [15] showed that the kernel estimators based density process corresponding to the general empirical process $\beta_n(x,1)$ converges weakly with the same normalization to the derivative of the limiting process in Theorem 1.1 of [20] that we quoted as Theorem A here.

We note that, with $F$ continuous, we have

$$\alpha_n(y,t) = \beta_n(Q(y),t), \qquad y, t \in [0,1],$$

and

$$\beta_n(x,t) = \alpha_n(F(x),t), \qquad x \in \mathbb{R}, t \in [0,1].$$

Hence, if $F$ is continuous, all strong and weak asymptotic results hold true simultaneously for both $\beta_n(x,t)$ and $\alpha_n(y,t)$.

For further reference, we spell out the weak convergence result that follows from Theorem A for $\alpha_n(y,t) = \beta_n(Q(y),t)$, $y, t \in [0,1]$, based on the induced sequence $\{U_n, n \geq 1\}$ as in (1.4).

COROLLARY A. *With $F$ continuous and $\tau$ and $D$ as in (1.3) and (1.5), respectively, as $n \to \infty$ we have*

$$\alpha_n(y,t) = \beta_n(Q(y),t) \xrightarrow{\mathcal{D}} \sqrt{\frac{2}{(2-\tau D)(1-\tau D)}} c_\tau(Q(y))Y_\tau(t)$$

$$= \sqrt{\frac{2}{(2-\tau D)(1-\tau D)}} J_\tau(y) Y_\tau(t)$$

in $D[0,1]^2$ that is equipped with the sup-norm, where, as before, $Y_\tau(t)$ is $1/\tau!$ times a Hermite process of rank $\tau$, given for each $t \in [0,1]$ as a multiple Wiener–Itô–Dobrushin integral as in, and right after, Theorem A.

In this paper we go further along these lines and establish strong approximations of the sequential uniform and general quantile processes, and of the sequential Bahadur–Kiefer processes as defined in (1.7) and (1.8) below. Moreover, we also study the sequential uniform Vervaat and Vervaat error processes of (1.10) and (1.11), respectively, along the same lines.

Since there is no simple relationship between $u_n(y,t)$ and $\gamma_n(y,t)$, following Csörgő and Révész [7] in the i.i.d. case along the lines of Csörgő and Szyszkowicz [11], here too we shall consider the normalized sequential general quantile process

$$(1.6) \quad \begin{aligned} \rho_n(y,t) &= f(Q(y))\gamma_n(y,t) = d_n^{-1}[nt]f(Q(y))(Q(y) - \widehat{Q_{[nt]}}(y)) \\ &= u_n(y,t) \frac{f(Q(y))}{f(Q(\theta_n(y,t)))}, \end{aligned}$$

where $0 \leq y, t \leq 1$, $|y - \theta_n(y,t)| \leq |y - \widehat{U_{[nt]}}(y)|$, provided that $F$ is an absolutely continuous distribution function with a strictly positive Lebesgue density function $f$ on the real line.

We define the stochastic processes

$$(1.7) \quad \begin{aligned} \{R_n^*(y,t), 0 &\leq y \leq 1, 0 \leq t \leq 1, n = 1, 2, \ldots\} \\ &= \{d_n(\alpha_n(y,t) - u_n(y,t)), 0 \leq y \leq 1, 0 \leq t \leq 1, n = 1, 2, \ldots\} \end{aligned}$$

and

$$(1.8) \quad \begin{aligned} \{R_n(y,t), 0 &\leq y \leq 1, 0 \leq t \leq 1, n = 1, 2, \ldots\} \\ &= \{d_n(\alpha_n(y,t) - \rho_n(y,t)), 0 \leq y \leq 1, 0 \leq t \leq 1, n = 1, 2, \ldots\} \\ &= \{d_n(\beta_n(Q(y),t) - \rho_n(y,t)), 0 \leq y \leq 1, 0 \leq t \leq 1, n = 1, 2, \ldots\}, \end{aligned}$$

which rhyme with the uniform and general Bahadur–Kiefer processes, respectively, in the i.i.d. case that enjoy some remarkable asymptotic properties (cf. [1, 23, 24] and Remark 1.2 below). For a review of, and contributions to, various aspects of this subject in the i.i.d. case, we refer to [4, 5, 6, 8, 9, 10, 11, 14, 30] and the references therein.

REMARK 1.2. It follows from the results of Kiefer [23, 24] that, in the i.i.d. case with $d_n = n^{1/2}$, $a_n R_n^*(\cdot, 1)$ cannot converge weakly in $D[0,1]$ to any nondegenerate random element of the latter space for any normalizing sequence $\{a_n\}$ of positive numbers. Vervaat [34, 35] argued this point in a crucially elegant way by showing that, in the i.i.d. case, in the space $C[0,1]$ (endowed with the uniform topology),

$$(1.9) \qquad V_n(s,1) := 2\int_0^s R_n^*(y,1)\,dy \xrightarrow{\mathcal{D}} B^2(s), \qquad n \to \infty,$$

where $B(\cdot)$ is a Brownian bridge. Accordingly then, if at all, $a_n R_n^*(\cdot,1)$ should converge weakly to a random element, say, $Y(\cdot)$, in $D[0,1]$, and we would then have to have the equality in distribution $\int_0^s Y(y)\,dy \stackrel{\mathcal{D}}{=} B^2(s)/2$, $0 \leq s \leq 1$. This, however, is impossible, for a Brownian bridge $B(\cdot)$ is almost surely nowhere differentiable. Vervaat [34, 35] established the above weak convergence of $V_n(\cdot,1)$ to the square of a Brownian bridge $B^2(\cdot)$ by showing that $\lim_{n\to\infty} \sup_{0\leq s\leq 1} |V_n(s,1) - \alpha_n^2(s,1)| = 0$ in probability. In view of this, one can think of the process $Q_n(s,1) := V_n(s,1) - \alpha_n^2(s,1)$, $0 \leq s \leq 1$, as the remainder term in the representation $V_n(s,1) = \alpha_n^2(s,1) + Q_n(s,1)$, $0 \leq s \leq 1$, of the uniform Vervaat process $V_n(\cdot,1)$ in terms of the square of the uniform empirical process $\alpha_n^2(\cdot,1)$. It is well known (cf., e.g., [36] for details and references) that $Q_n(\cdot,1)$ is asymptotically smaller than $\alpha_n^2(\cdot,1)$. Csörgő and Zitikis [12, 13, 14] and Csáki et al. [4] call this remainder term $Q_n(\cdot,1)$ the Vervaat error process, and study its strong and weak pointwise, sup-norm and $L_p$-norm asymptotic behavior for i.i.d. samples à la Kiefer [24] and Csörgő and Shi [9, 10]. Csörgő and Zitikis [13] and Csáki et al. [4] conclude that, just like the Bahadur–Kiefer process, in the i.i.d. case, $a_n Q_n(\cdot,1)$ cannot converge weakly to a nondegenerate random element in $D[0,1]$ for any sequence $\{a_n\}$ of positive real numbers.

In view of our discussion in Remark 1.2, based on $R_n^*(\cdot,\cdot)$ as in (1.7), we now introduce the integrated Bahadur–Kiefer process

$$(1.10) \quad V_n(s,t) = 2d_n^{-2}[nt]\int_0^s R_n^*(y,t)\,dy, \qquad 0 \leq s \leq 1, 0 \leq t \leq 1,$$

the so-called sequential uniform Vervaat process, and define the sequential Vervaat error process $Q_n(s,t)$ by

$$(1.11) \qquad Q_n(s,t) = V_n(s,t) - \alpha_n^2(s,t), \qquad 0 \leq s \leq 1, 0 \leq t \leq 1.$$

We shall see in this paper that, unlike in the i.i.d. case (cf. Remark 1.2), when appropriately normalized, the sequential Bahadur–Kiefer processes $R_n^*(\cdot,\cdot)$ and $R_n(\cdot,\cdot)$, as well as the sequential uniform Vervaat error process $Q_n(\cdot,\cdot)$, when based on long-range dependent sequences as in (1.2) and (1.4), do converge weakly in $D[0,1]^2$ (cf. Theorems 2.3, 3.2, 4.1 and 4.2), by first establishing strong approximations for these processes in sup-norm (cf. Theorems 2.1 and 2.2 and Propositions 3.1, 3.2 and 4.2). This new phenomenon in this context will be seen to be due to the limiting processes being Dehling–Taqqu type processes (cf. Remark 1.1), that is, multiplications of a nonrandom function by a random process which typically is a power of $Y_\tau(t)$, of Theorem A. Thus, via strong invariance, we arrive at functional limit theorems and laws of the iterated logarithm for the sequential Bahadur–Kiefer and the sequential uniform Vervaat error processes.

In Sections 2 and 3 we present strong invariance principles (approximations) for the sequential uniform Bahadur–Kiefer process and sequential uniform Vervaat error process of long-range dependent sequences as in (1.2) and (1.4), namely, for $R_n^*(y,t)$ and $Q_n(s,t)$ as in (1.7) and (1.11), respectively. Section 4 is devoted to establishing analogous statements for the sequential general Bahadur–Kiefer process $R_n(y,t)$ of (1.8) by examining the sup-norm distance between the sequential uniform quantile process $u_n(y,t)$ and the normalized sequential general quantile process $\rho_n(y,t)$ à la Csörgő and Révész [7] and Csörgő and Szyszkowicz [11]. The thus obtained results of Proposition 4.2 and Theorems 4.1 and 4.2 constitute a basis for studying quantiles, quantile and Bahadur–Kiefer processes in the context of long range dependent Gaussian subordinated processes.

The results obtained in this paper for long-range dependent sequences are analogs of those in the i.i.d. case in [4, 7, 8, 9, 10, 11, 12, 13].

REMARK 1.3. Further to long-range dependence, we note that long memory moving average models constitute an important and well-studied area of interest in time series analysis. In this regard, Koul and Surgailis [26] review various results on the asymptotic distribution of empirical processes of long memory moving averages with finite and infinite variance. Giraitis and Surgailis [22] discuss the uniform reduction principle for the empirical process of a long memory moving average process that generalizes the corresponding reduction principle of Dehling and Taqqu [20], which we also make fundamental use of in our Propositions 2.1 and 2.2. Thus, in principle, it should be possible to extend our present results to long memory moving average models as well. However, this extension is not within the immediate scope of the present paper. For a comprehensive study of empirical process techniques for dependent data in general, we refer to [17, 18].

## 2. Sequential uniform Bahadur–Kiefer process, strong approximations.

2.1. *Preliminaries.* Throughout this paper we assume that $\{X_n = G(\eta_n)\}$ and $\{U_n = F(G(\eta_n))\}$, $n \geq 1$, are as in (1.2) and (1.4), respectively, long-range dependent random sequences that are governed by the standard Gaussian random process $\{\eta_n\}$ which satisfies (1.1).

We first derive a strong approximation of the sequential general empirical process $\beta_n(x,t)$ by the process $c_\tau(x) \sum_{i=1}^{[nt]} H_\tau(\eta_i)/\tau!$, by changing the rate of convergence in Theorem 3.1 of [20] to fit our purposes in this exposition.

PROPOSITION 2.1. *Let $p$ be the smallest integer satisfying $\max(2, \tau, \frac{\tau D}{1-\tau D}) < p \leq \max(\frac{4-\tau D}{D}, \frac{4-\tau D}{1-\tau D})$. Then, as $n \to \infty$, we have*

$$\sup_{0 \leq t \leq 1} \sup_{-\infty < x < +\infty} \left| \beta_n(x, t) - d_n^{-1} c_\tau(x) \sum_{i=1}^{[nt]} H_\tau(\eta_i)/\tau! \right|$$

$$= O(n^{-\nu p/2 + \tau D/4 + \varepsilon}) \quad a.s.$$

*with any sufficiently small positive $\varepsilon$, where $\nu = \min(D, 1 - \tau D)/2$.*

PROOF. The proof is based on the well-known chaining argument of [20]. Hence, while studying the proof of Theorem 3.1 in [20], we shall only briefly indicate the extra steps that are needed for us to achieve our goal.

Let $S_n(k; x, y) = S_n(k; x) - S_n(k; y)$ $(-\infty < y \leq x < +\infty)$, where

$$S_n(k; x) = d_n^{-1} \sum_{i=1}^{k} \{I(X_i \leq x) - F(x) - c_\tau(x) H_\tau(\eta_i)/\tau!\}, \quad 1 \leq k \leq n.$$

Then we have

$$d_n S_n(k; x, y) = \sum_{i=1}^{k} \sum_{q=\tau+1}^{\infty} \frac{c_q(x) - c_q(y)}{q!} H_q(\eta_i).$$

This means that, for any fixed $-\infty < y \leq x < +\infty$, the Hermite rank of $d_n S_n(k; x, y)$ is at least $\tau + 1$.

First, for $\gamma(\cdot)$ as in (1.1), we assume that $\sup_{u \geq 1} |\gamma(u)| < \delta$, where $0 < \delta < (p-1)^{-1}$, and proceed as follows.

Via Proposition 4.2 of [32], one can verify that

$$\mathrm{E}|d_n S_n(k; x, y)|^p \leq C(p, \delta, x, y) \left\{ k \sum_{u=0}^{k} |\gamma(u)|^{\tau+1} \right\}^{p/2}$$

for some positive finite constant $C(p, \delta, x, y)$ depending on $p$, $\delta$ and $c_q(x) - c_q(y)$.

Since $|\frac{c_q(x) - c_q(y)}{q!}| \leq 2$ uniformly for any $-\infty < y \leq x < +\infty$, the proof of Lemma 4.5 and Proposition 4.2 of Taqqu [32] imply that the constant $C(p, \delta, x, y)$ must be a finite constant for any $-\infty < y \leq x < +\infty$. Letting $C(p, \delta, x, y) \leq C$, we get

$$\mathrm{E}|d_n S_n(k; x, y)|^p \leq C \left\{ k \sum_{u=0}^{k} |\gamma(u)|^{\tau+1} \right\}^{p/2}.$$

Suppose first $0 < D < (\tau+1)^{-1}$. Then, by (1.1), as $k \to \infty$,

$$k \sum_{u=0}^{k} |\gamma(u)|^{\tau+1} = O(k^{2-(\tau+1)D} L^{\tau+1}(k)).$$

When $D \geq (\tau+1)^{-1}$, $\sum_{u=0}^{k} |\gamma(u)|^{\tau+1}$ is slowly varying as $k \to \infty$ and, hence,

$$k \sum_{u=0}^{k} |\gamma(u)|^{\tau+1} = O(kL_0(k))$$

for some slowly varying function $L_0(\cdot)$ at infinity. Thus, we arrive at

$$(2.1) \quad \mathrm{E}|S_n(k;x,y)|^p \leq Ck^{-\nu p+\varepsilon}\left(\frac{d_k}{d_n}\right)^p \leq C\left(\frac{k}{n}\right)^{(1-\nu-\tau D/2)p} n^{-\nu p+\varepsilon},$$

with any sufficiently small positive $\varepsilon$ for any $-\infty < y \leq x < +\infty$, $1 \leq k \leq n$.

For any $s \geq 1$, define the partition as in [20],

$$-\infty = \pi_{0,s} < \pi_{1,s} < \cdots < \pi_{2^s,s} = +\infty$$

such that

$$\Lambda(\pi_{i,s}-) - \Lambda(\pi_{i-1,s}) \leq \Lambda(+\infty)2^{-s}, \qquad i = 1,\ldots,2^s,$$

$$\Lambda(x) = F(x) + \int_{\{G(y) \leq x\}} \frac{|H_\tau(y)|}{\tau!} \phi(y)\,dy,$$

where $\phi$ denotes the density function of the unit normal distribution.

Given $\zeta > 0$, let $K = [\log_2(C\zeta^{-1}nd_n^{-1})] + 1$. Next, for any $x \in \mathbb{R}$ and $s = 0, 1, \ldots, K$, define $j_s^x$ by

$$\pi_{j_s^x,s} \leq x < \pi_{j_s^x+1,s}.$$

One can then define a chain linking $-\infty$ to each point $x$ by

$$-\infty = \pi_{j_0^x,0} \leq \pi_{j_1^x,1} \leq \cdots \leq x < \pi_{j_K^x+1,K}.$$

Now using (2.1) instead of Lemma 3.1 of [20] and applying Chebyshev's inequality, along the same lines as those of the proof of Lemma 3.2 of [20], we obtain

$$\mathrm{P}\left\{\sup_{-\infty < x < +\infty} |S_n(k;x)| > \zeta\right\}$$

$$\leq \sum_{s=0}^{K} \mathrm{P}\left\{\sup_{-\infty < x < +\infty} |S_n(k;\pi_{j_s^x,s},\pi_{j_{s+1}^x,s+1})| > \zeta/(s+3)^2\right\}$$

$$+ \mathrm{P}\left\{d_n^{-1}\left|\sum_{i=1}^{n} H_\tau(\eta_i)\right| > 2^{K-1}\zeta/4\right\}$$

$$\leq C\left(\frac{k}{n}\right)^{(1-\nu-\tau D/2)p} n^{-\nu p+\varepsilon} \zeta^{-p} \sum_{s=0}^{K} 2^{s+1}(s+3)^{2p} + C\left(\frac{d_k}{d_n}\right)^p \zeta^{-p} 2^{-p(K-1)}$$

$$\leq C\left(\frac{k}{n}\right)^{(1-\nu-\tau D/2)p} n^{-\nu p+\varepsilon} \zeta^{-p} 2^K (K+3)^{2p+1}$$

$$+ C\left(\frac{k}{n}\right)^{(1-\tau D/2)p} n^{-\tau Dp/2+\varepsilon}$$

$$\leq Cn^{-\nu p+\tau D/2+\varepsilon}\left(\left(\frac{k}{n}\right)^{(1-\nu-\tau D/2)p}\zeta^{-p-\varepsilon} + \left(\frac{k}{n}\right)^{(1-\tau D/2)p}\right)$$

$$\leq Cn^{-\nu p+\tau D/2+\varepsilon}\left(\zeta^{-p-\varepsilon} + \left(\frac{k}{n}\right)^{(1-\tau D/2)p}\right)$$

for any $\zeta \in (0,1]$. The last inequality is due to the fact that $(1-\nu-\tau D/2)p > 1$.

Let $M_n(k) = \sup_{-\infty < x < +\infty} |S_n(k;x)|$, $M_n(k_1, k_2) = M_n(k_1) - M_n(k_2)$. On applying the above inequality, an appropriate variant of the proof of Theorem 3.1 of [20] leads to

$$P\left\{\max_{k\leq n} |M_n(k)| > \zeta\right\}$$

$$\leq \sum_{s=0}^{r} P\left\{\max_{j=1,\ldots,2^{r-s}} |M_n((j-1)2^s, j2^s)| > \frac{\zeta}{(s+2)^2}\right\}$$

$$\leq Cn^{-\nu p+\tau D/2+\varepsilon}\left\{\sum_{s=0}^{\log_2 n}(s+2)^{2p+2\varepsilon}\zeta^{-p-\varepsilon} + \sum_{s=0}^{\log_2 n} 2^{(s-r)(1-\tau D/2)p}\right\}$$

$$\leq Cn^{-\nu p+\tau D/2+\varepsilon_1}(1+\zeta^{-p-\varepsilon_1})$$

for some $\varepsilon_1 > \varepsilon$, where $r = \log_2 n$. This implies that, on assuming $\sup_{u\geq 1} |\gamma(u)| < \delta$ with $0 < \delta < (p-1)^{-1}$,

$$(2.2) \quad P\left\{\max_{k\leq n}\sup_{-\infty<x<+\infty} |S_n(k;x)| > \zeta\right\} \leq Cn^{-\nu p+\tau D/2+\varepsilon}(1+\zeta^{-p-\varepsilon}).$$

Now we proceed to establish (2.2) without the assumption $\sup_{u\geq 1}|\gamma(u)| < \delta$. Since $\gamma(u)$ tends to zero as $u \to \infty$, there exists a fixed integer $M = M(\delta) > 1$ such that $|\gamma(u)| < \delta$ for all $u \geq M$. Thus, without the assumption $\sup_{u\geq 1}|\gamma(u)| < \delta$, we merely have $|\gamma(u)| < \delta$ for $u \geq M$. Along the lines of the proof of Theorem 1 of [32], we obtain

$$\max_{1\leq k\leq n}\sup_{-\infty<x<\infty}|S_n(k;x)| = \max_{1\leq k\leq n}\sup_{-\infty<x<\infty}\left|d_n^{-1}\sum_{i=1}^{k}G^*(\eta_i;x)\right|$$

$$\leq \sum_{j=1}^{M}\max_{1\leq k\leq n}\sup_{-\infty<x<\infty}\left|d_n^{-1}\sum_{i=1}^{k'}G^*(\eta_{j+(i-1)M};x)\right|,$$

where $G^*(\eta_i;x) = I(X_i \leq x) - F(x) - c_\tau(x)H_\tau(\eta_i)/\tau!$.

Obviously, for each $j = 1, \ldots, M$, the correlations of the sequence $\{\eta_{j+(i-1)M}, i \geq 1\}$ are bounded by $\delta$ in absolute value. Let $S_n^*(k'; x) = d_n^{-1} \times \sum_{i=1}^{k'} G^*(\eta_{j+(i-1)M}; x)$. Therefore, (2.2) in our present context implies

$$P\left\{\max_{1 \leq k \leq n} \sup_{-\infty < x < +\infty} |S_n(k; x)| > \zeta\right\}$$

$$\leq \sum_{j=1}^{M} P\left\{\max_{1 \leq k \leq n} \sup_{-\infty < x < +\infty} |S_n^*(k'; x)| > \zeta\right\}$$

$$\leq C n^{-\nu p + \tau D/2 + \varepsilon}(1 + \zeta^{-p-\varepsilon}).$$

That is to say, (2.2) holds without the assumption $\sup_{u \geq 1} |\gamma(u)| < \delta$.

We now make use of (2.2) with $n = n_l = \min\{j : j \geq e^l\}$ and $\zeta = \zeta_l = \exp\{l(-\nu p/2 + \tau D/4 + \varepsilon)/(p + \varepsilon)\}$, $l = 0, 1, \ldots$. Then, by Borel–Cantelli lemma, there exists an integer $l_0$ such that, for any $l \geq l_0$,

$$\max_{k \leq n_l} \sup_{-\infty < x < +\infty} |S_{n_l}(k; x)| \leq \exp\{-l(\nu p/2 - \tau D/4 - \varepsilon)\} \quad \text{a.s.}$$

Let $n \geq e^{l_0}$ and let $l$ be the integer such that $n_{l-1} \leq n < n_l$. Since $e^{-l} \leq n^{-1}$ and $l \to \infty$ as $n \to \infty$, by definition of $d_n$ and that of a slowly varying function, we have

$$\sup_{-\infty < x < +\infty} |S_n(n; x)| \leq \frac{d_{n_l}}{d_n} \max_{k \leq n_l} \sup_{-\infty < x < +\infty} |S_{n_l}(k; x)| \leq C n^{-\nu p/2 + \tau D/4 + \varepsilon}.$$

This implies that, as $n \to \infty$,

$$\sup_{-\infty < x < +\infty} d_n^{-1} n^{\nu p/2 - \tau D/4 - \varepsilon} \left| d_n \beta_n(x, 1) - c_\tau(x) \sum_{i=1}^{n} H_\tau(\eta_i)/\tau! \right| = O(1) \quad \text{a.s.}$$

The latter, in turn, gives that, with fixed $t \in (0, 1]$ and $(nt) \to \infty$ as $n \to \infty$, we have

$$\sup_{-\infty < x < +\infty} \left| d_n \beta_n(x, t) - c_\tau(x) \sum_{i=1}^{[nt]} H_\tau(\eta_i)/\tau! \right|$$

$$= O(d_{[nt]}(nt)^{-\nu p/2 + \tau D/4 + \varepsilon})$$

$$= O((nt)^{1 - \nu p/2 - \tau D/4 + \varepsilon} L^{\tau/2}(nt)) \quad \text{a.s.,}$$

where the constant of $O(\cdot)$ is not a function of $t$ and, by our assumption for $p$, we see that the exponent of $(nt)^{1 - \nu p/2 - \tau D/4 + \varepsilon}$ is positive. Hence, without loss of generality, we can assume that the regularly varying function $(nt)^{1 - \nu p/2 - \tau D/4 + \varepsilon} L^{\tau/2}(nt)$ of positive exponent is a strictly monotone increasing regularly varying function of $(nt)$ (cf. 7 of Corollary 1.2.1 of [16] or Theorem 1.5.4 of [2]). Hence, on dividing both sides by $n^{1-\nu p/2-\tau D/4+\varepsilon}L^{\tau/2}(n)$, the right-hand side is seen to be a.s. bounded, independently of $t$. Consequently, we can take $\sup_{0 \leq t \leq 1}$ on the left-hand side, and thus arrive at the result of Proposition 2.1. $\square$

In the rest of this paper the marginal distribution function $F$ of $\{X_n\}$ in (1.2) is assumed to be continuous. We also assume the following:

ASSUMPTION A. $J_\tau(F(x))$ and the derivatives $J'_\tau(F(x))$, $J''_\tau(F(x))$ with $\tau$ as in (1.3) are uniformly bounded and
$$\sup_{0 < F(x) \leq \delta_n} |J_\tau(F(x))| = O(\delta_n)$$
for any sequence $\delta_n \to 0$ as $n \to \infty$.

REMARK 2.1. Since we assume that $F$ is continuous, if $J_\tau(0) = 0$, it follows that
$$\sup_{0 < F(x) \leq \delta_n} |J_\tau(F(x))| = \sup_{0 < F(x) \leq \delta_n} |J_\tau(F(x)) - J_\tau(0)|$$
$$\leq \sup_{0 < F(x) \leq \delta_n} F(x) \cdot \sup_{0 < F(\theta) \leq \delta_n} |J'_\tau(F(\theta))| = O(\delta_n).$$

Moreover, if we take $G$ as $G = F^{-1}\Phi$, we see that $J_1(F(x)) = -\phi(\Phi^{-1}(F(x))) \neq 0$ for any $F(x) \in (0,1)$, where $\phi$, $\Phi^{-1}$ denote, respectively, the density function and the quantile function of the unit normal distribution function $\Phi$. This means that in this case $\tau = 1$, and elementary calculations show that Assumption A holds automatically. Specifically, let $G(x) = x$, and this is a special case of $G = F^{-1}\Phi$, since now $F = \Phi$. Then, for this function $G$, $\tau = 1$ and $J_1(\Phi(x))$ satisfies Assumption A.

For the sake of first approximating the sequential uniform empirical and quantile processes $\alpha_n(y,t)$ and $u_n(y,t)$, we define the two-time parameter stochastic process $\{V(y,nt); 0 \leq y \leq 1, 0 \leq t \leq 1, n \geq 1\}$ by

$$(2.3) \qquad V(y,nt) = J_\tau(y) \sum_{i=1}^{[nt]} H_\tau(\eta_i)/\tau!,$$

and, as an immediate consequence of Proposition 2.1, we conclude the following strong approximation for the sequential uniform empirical process $\alpha_n(y,t)$.

COROLLARY 2.1. *Under the assumptions of Proposition 2.1, we have*
$$\sup_{0 \leq t \leq 1} \sup_{0 \leq y \leq 1} |\alpha_n(y,t) - d_n^{-1}V(y,nt)| = O(n^{-\nu p/2 + \tau D/4 + \varepsilon}) \qquad a.s.$$
*with any sufficiently small positive $\varepsilon$, where $\nu = \min(D, 1-\tau D)/2$.*

Let $\kappa_{1\tau} = \sup_{0 \leq y \leq 1} |J_\tau(y)|$, $\kappa_{2\tau} = \sup_{0 \leq y \leq 1} |J_\tau(y) \cdot J'_\tau(y)|$, $\kappa_{3\tau} = \sup_{0 \leq y \leq 1} |J^2_\tau(y) \cdot J'_\tau(y)|$. Via Assumption A, we conclude $0 < \kappa_{1\tau}, \kappa_{2\tau}, \kappa_{3\tau} < \infty$. Moreover, if we take $G = F^{-1}\Phi$, by Remark 2.1, it is easy to check that $\kappa_{11} = 1/(2\pi)^{1/2}$, $\kappa_{21} = 1/(2\pi e)^{1/2}$ and $\kappa_{31} = 1/\{2\pi(2e)^{1/2}\}$.

The process $V(y, nt)$ defined in (2.3) that is approximating $\alpha_n(y, t)$ as in Corollary 2.1 can also be used to approximate the sequential uniform quantile process $u_n(y, t)$. Namely, we have the following:

PROPOSITION 2.2.  *Let $p$ be the smallest integer satisfying $\max(3\tau, \frac{3\tau D}{1-\tau D}) < p \leq \max(\frac{4-\tau D}{D}, \frac{4-\tau D}{1-\tau D})$. Suppose Assumption A holds. Then under the assumptions of Corollary 2.1, as $n \to \infty$, we have*

$$\sup_{0 \leq t \leq 1} \sup_{0 \leq y \leq 1} |u_n(y, t) - d_n^{-1} V(y, nt)|$$
(2.4)
$$= O(n^{-\tau D/2} L^{\tau/2}(n)(\log \log n)^\tau) \quad a.s.$$

PROOF.  Note that

$$u_n(y, t) = d_n^{-1}[nt]\{\widehat{E_{[nt]}}(\widehat{U_{[nt]}}(y)) - \widehat{U_{[nt]}}(y)\} - d_n^{-1}[nt]\{\widehat{E_{[nt]}}(\widehat{U_{[nt]}}(y)) - y\}$$
$$= \alpha_n(\widehat{U_{[nt]}}(y), t) - d_n^{-1}[nt]\{\widehat{E_{[nt]}}(\widehat{U_{[nt]}}(y)) - y\},$$

and it is easy to see that

$$0 \leq \sup_{0 \leq y \leq 1} |\widehat{E_n}(\widehat{U_{[nt]}}(y)) - y| \leq 1/[nt].$$

Thus, we have

$$\sup_{0 \leq t \leq 1} \sup_{0 \leq y \leq 1} |u_n(y, t) - \alpha_n(y, t)|$$
(2.5)
$$= \sup_{0 \leq t \leq 1} \sup_{0 \leq y \leq 1} |\alpha_n(\widehat{U_{[nt]}}(y), t) - \alpha_n(y, t)| + O(d_n^{-1}).$$

Applying Corollary 2.1, estimating the right-hand side of (2.5), we obtain

$$\sup_{0 \leq t \leq 1} \sup_{0 \leq y \leq 1} |\alpha_n(\widehat{U_{[nt]}}(y), t) - \alpha_n(y, t)|$$
(2.6)
$$= \sup_{0 \leq t \leq 1} \sup_{0 \leq y \leq 1} d_n^{-1} |V(\widehat{U_{[nt]}}(y), nt) - V(y, nt)|$$
$$+ O(n^{-\nu p/2 + \tau D/4 + \varepsilon}) \quad a.s.$$

Hence, we need to study the size of the random increments of the process $V(y, nt)$.

The Mori–Oodaira LIL [28] yields

$$\limsup_{n\to\infty} n^{\tau D/2-1}(L(n)\log\log n)^{-\tau/2} \sup_{0\leq t\leq 1}\left|\sum_{i=1}^{[nt]} H_\tau(\eta_i)/\tau!\right| \quad (2.7)$$
$$= \frac{2^{(\tau+1)/2}}{\sqrt{\tau!(2-\tau D)(1-\tau D)}} \quad \text{a.s.}$$

Hence, by (2.3) and the fact that $0 < \kappa_{1\tau} < \infty$, we have

$$\limsup_{n\to\infty}(\log\log n)^{-\tau/2} \sup_{0\leq t\leq 1}\sup_{0\leq y\leq 1} d_n^{-1}|V(y,nt)| \quad (2.8)$$
$$= \frac{2^{(\tau+1)/2}\kappa_{1\tau}}{\sqrt{\tau!(2-\tau D)(1-\tau D)}} \quad \text{a.s.}$$

Consequently, via Corollary 2.1, we conclude

$$\limsup_{n\to\infty}(\log\log n)^{-\tau/2} \sup_{0\leq t\leq 1}\sup_{0\leq y\leq 1}|\alpha_n(y,t)| = \frac{2^{(\tau+1)/2}\kappa_{1\tau}}{\sqrt{\tau!(2-\tau D)(1-\tau D)}} \quad \text{a.s.,}$$

and this in turn gives

$$\limsup_{n\to\infty}(\log\log n)^{-\tau/2} \sup_{0\leq t\leq 1}\sup_{0\leq y\leq 1}|u_n(y,t)| \quad (2.9)$$
$$= \frac{2^{(\tau+1)/2}\kappa_{1\tau}}{\sqrt{\tau!(2-\tau D)(1-\tau D)}} \quad \text{a.s.}$$

on account of

$$\sup_{0\leq t\leq 1}\sup_{0\leq y\leq 1}|\alpha_n(y,t)| = \sup_{0\leq t\leq 1}\sup_{0\leq y\leq 1}|u_n(y,t)|.$$

On the other hand, by the mean value theorem, we arrive at

$$|J_\tau(\widehat{U}_n(y)) - J_\tau(y)| = |\widehat{U}_n(y) - y||J'_\tau(\theta_{1n}(y))|,$$

where $|y - \theta_{1n}(y)| \leq |\widehat{U}_n(y) - y|$. Now (2.9) with $t=1$ implies that, as $n\to\infty$,

$$\sup_{0\leq y\leq 1}|\widehat{U}_n(y) - y| = \sup_{0\leq y\leq 1} d_n n^{-1}|u_n(y,1)| \quad (2.10)$$
$$= O((n^{-D}L(n)\log\log n)^{\tau/2}) \to 0$$

almost surely [we note in passing that (2.10) is just a *Glivenko–Cantelli theorem* with rates of convergence in terms of the long-range dependent sequence as in (1.4)]. Thus, by Assumption A, as $n\to\infty$, we arrive at

$$\sup_{0\leq y\leq 1}|J_\tau(\widehat{U}_n(y)) - J_\tau(y)|$$
$$= O((n^{-D}L(n)\log\log n)^{\tau/2}) \quad \text{a.s.}$$

The latter combined with (2.7) for $t = 1$ yields

$$\sup_{0 \leq y \leq 1} d_n^{-1} |V(\widehat{U}_n(y), n) - V(y, n)|$$

(2.11)
$$= O(n^{-\tau D/2} L^{\tau/2}(n) (\log \log n)^\tau) \quad \text{a.s.}$$

Using (2.5)–(2.6), (2.11) and our assumption for $p$, we arrive at

$$\sup_{0 \leq y \leq 1} |u_n(y, 1) - \alpha_n(y, 1)|$$

(2.12)
$$= O(n^{-\tau D/2} L^{\tau/2}(n) (\log \log n)^\tau) \quad \text{a.s.}$$

Now (2.12) combined with Corollary 2.1 with $t = 1$ yields

$$\sup_{0 \leq y \leq 1} |u_n(y, 1) - d_n^{-1} V(y, n)|$$

$$= O(n^{-\tau D/2} L^{\tau/2}(n) (\log \log n)^\tau) \quad \text{a.s.}$$

On multiplying through by $d_n$ and then applying a similar argument as used at the end of the proof of Proposition 2.1, we conclude (2.4). □

Next, in view of (2.5) and (2.6), we establish the exact size of the random increments of the process $V(y, nt)$ for convenient use later on.

PROPOSITION 2.3. *Under the assumptions of Proposition* 2.2, *we have*

$$\limsup_{n \to \infty} n^{\tau D - 1} (L(n) \log \log n)^{-\tau} \sup_{0 \leq t \leq 1} \sup_{0 \leq y \leq 1} |V(\widehat{U}_{[nt]}(y), nt) - V(y, nt)|$$

$$= \frac{2^{\tau+1} \kappa_{2\tau}}{\tau!(2 - \tau D)(1 - \tau D)} \quad a.s.$$

PROOF. We note that

$$J_\tau(\widehat{U}_n(y)) - J_\tau(y)$$
$$= J_\tau'(y)(\widehat{U}_n(y) - y) + \tfrac{1}{2}(\widehat{U}_n(y) - y)^2 J_\tau''(\theta_{2n}(y))$$
$$= -J_\tau'(y) n^{-1} V(y, n) + J_\tau'(y) d_n n^{-1} (d_n^{-1} V(y, n) - u_n(y))$$
$$+ \tfrac{1}{2}(\widehat{U}_n(y) - y)^2 J_\tau''(\theta_{2n}(y)),$$

where $|y - \theta_{2n}(y)| \leq |\widehat{U}_n(y) - y|$. Consequently, by (2.4) with $t = 1$ and (2.10), we obtain

$$\sup_{0 \leq y \leq 1} |J_\tau'(y) d_n n^{-1} (u_n(y) - d_n^{-1} V(y, n))|$$

$$= O((n^{-D} L(n) \log \log n)^\tau) \quad \text{a.s.}$$

and
$$\sup_{0\leq y\leq 1} |\tfrac{1}{2}(\widehat{U_n}(y)-y)^2 J_\tau''(\theta_{2n}(y))|$$
$$= O((n^{-D}L(n)\log\log n)^\tau) \quad \text{a.s.}$$

Hence,
$$\sup_{0\leq y\leq 1} |J_\tau(\widehat{U_n}(y)) - J_\tau(y) + J_\tau'(y)n^{-1}V(y,n)|$$
$$= O((n^{-D}L(n)\log\log n)^\tau) \quad \text{a.s.}$$

Now (2.7) with $t=1$ implies
$$\limsup_{n\to\infty} (n^{-D}L(n)\log\log n)^{-\tau/2} \sup_{0\leq y\leq 1} |J_\tau(\widehat{U_n}(y)) - J_\tau(y)|$$
$$= \frac{2^{(\tau+1)/2}\kappa_{2\tau}}{\sqrt{\tau!(2-\tau D)(1-\tau D)}}$$

almost surely and, again by (2.7), we conclude that
$$\limsup_{n\to\infty} n^{\tau D-1}(L(n)\log\log n)^{-\tau} \sup_{0\leq y\leq 1} |V(\widehat{U_n}(y),n) - V(y,n)|$$
$$= \frac{2^{\tau+1}\kappa_{2\tau}}{\tau!(2-\tau D)(1-\tau D)} \quad \text{a.s.}$$

Hence, we have, with $t\in(0,1)$ fixed, as $(nt)\to\infty$,
$$\sup_{0\leq y\leq 1} |V(\widehat{U_{[nt]}}(y),nt) - V(y,nt)|$$
$$= \left(\frac{2^{\tau+1}\kappa_{2\tau}}{\tau!(2-\tau D)(1-\tau D)} + o(1)\right)(nt)^{1-\tau D}L^\tau(nt)(\log\log(nt))^\tau \quad \text{a.s.,}$$

and hence, on dividing both sides by $n^{1-\tau D}L^\tau(n)(\log\log n)^\tau$ and assuming without loss of generality that the regularly varying function $n^{1-\tau D}L^\tau(n)$ of positive exponent is strictly monotone increasing, the right-hand side is seen to be a.s. bounded and independent of $t$. Thus, on taking $\sup_{0\leq t\leq 1}$ on the left-hand side, we conclude the proof of Proposition 2.3. □

PROPOSITION 2.4. *Under the assumptions of Proposition* 2.2, *as* $n\to\infty$, *we have*
$$\sup_{0\leq t\leq 1}\sup_{0\leq y\leq 1} |V(\widehat{U_{[nt]}}(y),nt) - V(y - [nt]^{-1}V(y,nt),nt)|$$
$$= O(n^{1-3\tau D/2}(L(n)\log\log n)^{3\tau/2}) \quad \text{a.s.}$$

*or, equivalently,*

$$\sup_{0\leq t\leq 1}\sup_{0\leq y\leq 1}|(V(\widehat{U_{[nt]}}(y),nt)-V(y,nt))$$

$$-(V(y-[nt]^{-1}V(y,nt),nt)-V(y,nt))|$$

$$=O(n^{1-3\tau D/2}(L(n)\log\log n)^{3\tau/2}) \quad a.s.$$

PROOF. Notice that

$$V(\widehat{U_n}(y),n)=V(y-n^{-1}V(y,n)-\Delta_n(y),n),$$

where $\Delta_n(y)=d_n n^{-1}(u_n(y,1)-d_n^{-1}V(y,n))$. By Proposition 2.2 with $t=1$, we get

$$\sup_{0\leq y\leq 1}|\Delta_n(y)|=O((n^{-D}L(n)\log\log n)^\tau) \quad a.s.$$

Consequently, along the lines of the proof for (2.11), we obtain

$$\sup_{0\leq y\leq 1}|V(\widehat{U_n}(y),n)-V(y-n^{-1}V(y,n),n)|$$

$$=O(n^{1-3\tau D/2}L^{3\tau/2}(n)(\log\log n)^{3\tau/2}) \quad a.s.$$

This also completes the proof of Proposition 2.4 by using a similar argument as at the end of the proof of Proposition 2.3. □

PROPOSITION 2.5. *Under the assumptions of Proposition 2.2, we have*

$$\sup_{0\leq t\leq 1}\sup_{0\leq y\leq 1}\left|V(y-[nt]^{-1}V(y,nt),nt)\right.$$

(2.13)
$$\left.-V(y,nt)+[nt]^{-1}V(y,nt)J'_\tau(y)\sum_{i=1}^{[nt]}H_\tau(\eta_i)/\tau!\right|$$

$$=O(n^{1-3\tau D/2}(L(n)\log\log n)^{3\tau/2}) \quad a.s., \ n\to\infty,$$

*and*

$$\limsup_{n\to\infty} n^{\tau D-1}(L(n)\log\log n)^{-\tau}$$

$$\times\sup_{0\leq t\leq 1}\sup_{0\leq y\leq 1}|V(y-[nt]^{-1}V(y,nt),nt)-V(y,nt)|$$

(2.14)
$$=\limsup_{n\to\infty} n^{\tau D-1}(L(n)\log\log n)^{-\tau}$$

$$\times\sup_{0\leq t\leq 1}\sup_{0\leq y\leq 1}\left|[nt]^{-1}V(y,nt)J'_\tau(y)\sum_{i=1}^{[nt]}H_\tau(\eta_i)/\tau!\right|$$

$$=\frac{2^{\tau+1}\kappa_{2\tau}}{\tau!(2-\tau D)(1-\tau D)} \quad a.s.$$

BAHADUR–KIEFER AND VERVAAT ERROR PROCESSES 19

PROOF. By (2.8) and (2.10), respectively, as $n \to \infty$, we have

$$\sup_{0 \leq y \leq 1} n^{-1}|V(y,n)| = O((n^{-D}L(n)\log\log n)^{\tau/2}) \qquad \text{a.s.}$$

and

$$\sup_{0 \leq y \leq 1} |\widehat{U_n}(y) - y| = O((n^{-D}L(n)\log\log n)^{\tau/2}) \qquad \text{a.s.}$$

Hence, along the lines of the proof of Proposition 2.3, we first obtain (2.13) and (2.14) with $t = 1$, and then a similar argument as at the end of the proof of Proposition 2.3 yields (2.13) and (2.14) as stated. $\square$

2.2. *Strong approximations of sequential uniform Bahadur–Kiefer processes.* A direct application of Corollary 2.1 and (2.5) leads to a strong approximation for the sequential uniform Bahadur–Kiefer process $R_n^*(y,t)$.

THEOREM 2.1. *Under the assumptions of Corollary* 2.1, *as $n \to \infty$, we have*

$$\sup_{0 \leq t \leq 1} \sup_{0 \leq y \leq 1} |R_n^*(y,t) - (V(y,nt) - V(\widehat{U_{[nt]}}(y), nt))|$$

$$= \sup_{0 \leq t \leq 1} \sup_{0 \leq y \leq 1} |d_n(\alpha_n(y,t) - u_n(y,t)) - (V(y,nt) - V(\widehat{U_{[nt]}}(y), nt))|$$

$$= O(n^{1-\nu p/2 - \tau D/4 + \varepsilon} L^{\tau/2}(n)) \qquad \text{a.s.}$$

Next we reformulate Theorem 2.1 as follows.

THEOREM 2.2. *Under the assumptions of Proposition* 2.2, *as $n \to \infty$, we have*

$$\sup_{0 \leq t \leq 1} \sup_{0 \leq y \leq 1} \left| [nt] R_n^*(y,t) - J_\tau(y) J'_\tau(y) \left( \sum_{i=1}^{[nt]} H_\tau(\eta_i)/\tau! \right)^2 \right|$$

$$= O(n^{2-\nu p/2 - \tau D/4 + \varepsilon} L^{\tau/2}(n)) \qquad \text{a.s.}$$

PROOF. Propositions 2.4, 2.5 and Theorem 2.1 imply the result. $\square$

These strong approximations readily yield weak convergence and laws of the iterated logarithm for the process $R_n^*(y,t)$.

THEOREM 2.3. *Under the assumptions of Proposition* 2.2, *as $n \to \infty$, we have*

$$n^{\tau D - 2} L^{-\tau}(n)[nt] R_n^*(y,t) \xrightarrow{\mathcal{D}} \frac{2}{(2-\tau D)(1-\tau D)} J_\tau(y) J'_\tau(y) Y_\tau^2(t)$$

in the space $D[0,1]^2$, equipped with the sup-norm, where $Y_\tau(t)$ is as in Theorem A.

PROOF. From Theorem 5.6 of [33], as $n \to \infty$, we conclude

$$d_n^{-1} \sum_{i=1}^{[nt]} H_\tau(\eta_i)/\tau! \xrightarrow{\mathcal{D}} \sqrt{\frac{2}{(2-\tau D)(1-\tau D)}} Y_\tau(t)$$

in $D[0,1]$. Now Theorem 2.3 follows from Theorem 2.2.　□

In light of Theorems 2.2 and 2.3, we have the following:

THEOREM 2.4.　*Under the assumptions of Proposition 2.2, we have*

(2.15)
$$\limsup_{n \to \infty} n^{\tau D - 2}(L(n) \log \log n)^{-\tau} \sup_{0 \le t \le 1} \sup_{0 \le y \le 1} |[nt]R_n^*(y,t)|$$
$$= \frac{2^{\tau+1}\kappa_{2\tau}}{\tau!(2-\tau D)(1-\tau D)} \quad a.s.,$$

*as well as*

(2.16)
$$n^{\tau D - 2} L^{-\tau}(n) \sup_{0 \le t \le 1} \sup_{0 \le y \le 1} |[nt]R_n^*(y,t)|$$
$$\xrightarrow{\mathcal{D}} \frac{2\kappa_{2\tau}}{(2-\tau D)(1-\tau D)} \sup_{0 \le t \le 1} Y_\tau^2(t), \qquad n \to \infty.$$

PROOF. Equation (2.15) follows from Theorem 2.2 and the law of the iterated logarithm (2.7) for $\sum_{i=1}^{[nt]} H_\tau(\eta_i)/\tau!$. As to (2.16), it results from Theorem 2.3 directly.　□

Denote the $L_p$-norm of a function $f$ on $[0,1]^2$ by

$$\|f\|_p = \left( \int_0^1 \int_0^1 |f(y,t)|^p \, dy \, dt \right)^{1/p}, \qquad 1 \le p < \infty.$$

A straightforward $L_p$-version of Theorem 2.2 for the sequential uniform Bahadur–Kiefer process $R_n^*(y,t)$ results in the following:

THEOREM 2.5.　*Under the assumptions of Proposition 2.2, we have*

$$\limsup_{n \to \infty} n^{\tau D - 2}(L(n) \log \log n)^{-\tau} \|[nt]R_n^*\|_p = \frac{2^{\tau+1}\|J_\tau J_\tau'\|_p}{\tau!(2-\tau D)(1-\tau D)} \qquad a.s.,$$

*as well as*

$$n^{\tau D - 2} L^{-\tau}(n) \|[nt]R_n^*\|_p \xrightarrow{\mathcal{D}} \frac{2\|J_\tau J_\tau'\|_p}{(2-\tau D)(1-\tau D)} \|Y_\tau^2\|_p, \qquad n \to \infty.$$

This is in contrast to the $L_p$-theory of the Bahadur–Kiefer process in the i.i.d. case in [9, 10], which deviates substantially from its [23, 24] sup-norm theory. For a review of this matter, we refer to [4]. For the sake of comparison to the latter theories, Theorems 2.1–2.5 above should be read with $t = 1$. For strong approximations in sup-norm of the sequential uniform Bahadur–Kiefer process in the i.i.d. case, we refer to [11].

**3. Asymptotics of the uniform Vervaat error process.** In support of studying the sequential uniform Vervaat error process, we first derive the weak convergence of the sequential uniform Vervaat process $V_n(\cdot,\cdot)$ [cf. (1.10)]. This can be easily done via Theorems 2.2 and 2.3.

THEOREM 3.1. *Under the assumptions of Proposition 2.2, as $n \to \infty$, we have*

$$V_n(s,t) \xrightarrow{\mathcal{D}} \frac{2}{(2-\tau D)(1-\tau D)} J_\tau^2(s) Y_\tau^2(t)$$

*in the space $D[0,1]^2$, equipped with the sup-norm, where $Y_\tau(t)$ is as in Theorem A.*

PROOF. Theorem 2.3 and integration by parts yield

$$\begin{aligned}
V_n(s,t) &= 2d_n^{-2}[nt] \int_0^s R_n^*(y,t)\, dy \\
&\xrightarrow{\mathcal{D}} \frac{4}{(2-\tau D)(1-\tau D)} \left( \int_0^s J_\tau(y) J_\tau'(y)\, dy \right) Y_\tau^2(t) \\
&= \frac{2}{(2-\tau D)(1-\tau D)} J_\tau^2(s) Y_\tau^2(t). \qquad \square
\end{aligned}$$

Theorem 3.1 and Corollary A imply that the sequential uniform Vervaat process $V_n(s,t)$ and the process $\alpha_n^2(s,t)$ have the same weak limiting process. Thus, just as in the i.i.d. case, it makes sense to consider the deviation of the two processes, that is, the sequential uniform Vervaat error process $Q_n$ as in (1.11). Unlike in the i.i.d. case (cf. [4]), we shall see that $Q_n(s,1)$, as well as its sequential version $Q_n(s,t)$, do converge weakly and, in particular, to a random process which is a multiplication of a nonrandom function by the cube of the random process $Y_\tau(t)$ of Theorem A.

PROPOSITION 3.1. *Under the assumptions of Proposition 2.2, we have*

$$\sup_{0 \leq t \leq 1} \sup_{0 \leq s \leq 1} [nt]|Q_n(s,t) - Z_n(s,t)| = O(n^{1-\nu p/2 + \tau D/4 + \varepsilon} (\log \log n)^{\tau/2}) \qquad a.s.,$$

*where* $\{Z_n(s,t), 0 \leq s,t \leq 1, n = 1,2,\ldots\}$ *is defined by*

$$Z_n(s,t) = 2d_n^{-2} V(s, nt)$$
(3.1)
$$\times \int_0^1 (V(s - w[nt]^{-1}V(s,nt), nt) - V(s,nt))\, dw.$$

PROOF. We proceed à la the lines of the proofs of Lemmas 3.1 and 3.2 of [4]. Let

$$A_n(s,t) = 2d_n^{-1}[nt] \int_{\widehat{U_{[nt]}}(s)}^s (\alpha_n(y,t) - \alpha_n(s,t))\, dy,$$
(3.2)
$$0 \leq s,t \leq 1, n = 1,2,\ldots.$$

It follows from Lemma 3.1 of [4] that

$$Q_n(s,t) = A_n(s,t) - d_n^{-2}(R_n^*(s,t))^2.$$

Now (2.15) with $t = 1$ yields that, when $n \to \infty$,

$$\sup_{0 \leq s \leq 1} n|Q_n(s,1) - A_n(s,1)| = O(n^{1-\tau D} L^\tau(n)(\log\log n)^{2\tau}) \quad \text{a.s.}$$

In similar fashion as at the end of the proof of Proposition 2.1, as $n \to \infty$ we get

$$\sup_{0 \leq t \leq 1} \sup_{0 \leq s \leq 1} [nt]|Q_n(s,t) - A_n(s,t)| = O(n^{1-\tau D} L^\tau(n)(\log\log n)^{2\tau}) \quad \text{a.s.}$$

Hence, it suffices to show that, as $n \to \infty$,

$$\sup_{0 \leq t \leq 1} \sup_{0 \leq s \leq 1} [nt]|A_n(s,t) - Z_n(s,t)| = O(n^{1-\nu p/2 + \tau D/4 + \varepsilon}(\log\log n)^{\tau/2}) \quad \text{a.s.}$$

Changing variable $y = s - w(s - \widehat{U_{[nt]}}(s)) = s - w[nt]^{-1} d_n u_n(s,t)$ in (3.2), we get

$$A_n(s,t) = 2u_n(s,t) \int_0^1 (\alpha_n(s - w[nt]^{-1} d_n u_n(s,t)) - \alpha_n(s,t))\, dw.$$

Corollary 2.1 and (2.9), as $n \to \infty$, yield

$$A_n(s,t) = 2d_n^{-1} u_n(s,t)$$
(3.3)
$$\times \int_0^1 (V(s - w[nt]^{-1} d_n u_n(s,t), nt) - V(s,nt))\, dw$$
$$+ O(n^{-\nu p/2 + \tau D/4 + \varepsilon}(\log\log n)^{\tau/2}) \quad \text{a.s.},$$

uniformly in $s,t \in [0,1]$. For all $0 \leq w \leq 1$, according to Proposition 2.4, as $n \to \infty$ we have, uniformly in $s,t \in [0,1]$,

$$V(s - w[nt]^{-1} d_n u_n(s,t), nt) = V(s - w[nt]^{-1} V(s,nt), nt)$$
$$+ O(n^{1-3\tau D/2}(L(n)\log\log n)^{3\tau/2}) \quad \text{a.s.}$$

Inserting this into (3.3) and applying (2.9) again, we obtain, uniformly in $s, t \in [0, 1]$,

$$A_n(s,t) = 2d_n^{-1} u_n(s,t) \int_0^1 (V(s - w[nt]^{-1} V(s, nt), nt) - V(s, nt))\, dw$$
$$+ O(n^{-\nu p/2 + \tau D/4 + \varepsilon}(\log \log n)^{\tau/2})$$
$$+ O(n^{-\tau D} L^\tau(n)(\log \log n)^{2\tau}) \quad \text{a.s.}$$

Consequently, as $n \to \infty$, uniformly in $s, t \in [0, 1]$,

(3.4)
$$A_n(s,t) = 2d_n^{-1} u_n(s,t)$$
$$\times \int_0^1 (V(s - w[nt]^{-1} V(s, nt), nt) - V(s, nt))\, dw$$
$$+ O(n^{-\nu p/2 + \tau D/4 + \varepsilon}(\log \log n)^{\tau/2}) \quad \text{a.s.}$$

Now, from Proposition 2.2, as $n \to \infty$,

(3.5) $\quad 2d_n^{-1} u_n(s,t) = 2d_n^{-2} V(s, nt) + O(n^{-1}(\log \log n)^\tau) \quad \text{a.s.}$

uniformly in $0 \le s, t \le 1$. On the other hand, applying (2.14) to the integrand in (3.4), we arrive at

$$\int_0^1 (V(s - w[nt]^{-1} V(s, nt), nt) - V(s, nt))\, dw$$
$$= O(n^{1 - \tau D}(L(n) \log \log n)^\tau) \quad \text{a.s.}$$

uniformly in $0 \le s, t \le 1$. Inserting this and (3.5) into (3.4) yields that, as $n \to \infty$,

$$[nt]|A_n(s,t) - Z_n(s,t)| = O(n^{1 - \nu p/2 + \tau D/4 + \varepsilon}(\log \log n)^{\tau/2}) \quad \text{a.s.}$$

uniformly in $0 \le s, t \le 1$. This concludes the proof of Proposition 3.1. □

Due to Proposition 2.5, we present the following conclusion.

PROPOSITION 3.2. *Under the assumptions of Proposition 2.2, we have*

$$\sup_{0 \le t \le 1} \sup_{0 \le s \le 1} \left| Z_n(s,t) + 2d_n^{-2}[nt]^{-1}(V(s,nt))^2 J'_\tau(s) \sum_{i=1}^{[nt]} H_\tau(\eta_i)/\tau! \right|$$
$$= O(n^{-\tau D} L^\tau(n)(\log \log n)^{2\tau}) \quad a.s.$$

PROOF. By (3.1) and Proposition 2.5, we obtain

$$\sup_{0\le t\le 1}\sup_{0\le s\le 1}\left|Z_n(s,t)+2d_n^{-2}[nt]^{-1}(V(s,nt))^2 J'_\tau(s)\sum_{i=1}^{[nt]}H_\tau(\eta_i)/\tau!\right|$$

$$=\sup_{0\le t\le 1}\sup_{0\le s\le 1}\left|2d_n^{-2}V(s,nt)\right.$$

$$\times\int_0^1\left\{V(s-w[nt]^{-1}V(s,nt),nt)-V(s,nt)\right.$$

$$\left.\left.+w[nt]^{-1}V(s,nt)J'_\tau(s)\sum_{i=1}^{[nt]}H_\tau(\eta_i)/\tau!\right\}dw\right|$$

$$\le\sup_{0\le t\le 1}\sup_{0\le s\le 1}|2d_n^{-2}V(s,nt)|$$

$$\times\sup_{0\le t\le 1}\sup_{0\le s,w\le 1}\left|V(s-w[nt]^{-1}V(s,nt),nt)\right.$$

$$\left.-V(s,nt)+w[nt]^{-1}V(s,nt)J'_\tau(s)\sum_{i=1}^{[nt]}H_\tau(\eta_i)/\tau!\right|$$

$$=O(n^{-\tau D}L^\tau(n)(\log\log n)^{2\tau})\qquad\text{a.s.}$$

This completes the proof. □

The main conclusions of this section are as follows.

THEOREM 3.2. *Under the assumptions of Proposition 2.2, as $n\to\infty$, we have*
$$n^{\tau D/2-1}L^{-\tau/2}(n)[nt]Q_n(s,t)\xrightarrow{\mathcal{D}}2^{5/2}((2-\tau D)(1-\tau D))^{-3/2}J_\tau^2(s)J'_\tau(s)Y_\tau^3(t)$$
*in the space $D[0,1]^2$, equipped with the sup-norm, where $Y_\tau(t)$ is as in Theorem A.*

PROOF. It follows from Theorem 5.6 of [33] and Propositions 3.1 and 3.2. □

As a consequence of Propositions 3.1 and 3.2 and Theorem 3.2, we have the following results.

THEOREM 3.3. *Under the conditions of Proposition 2.2, we have*
$$\limsup_{n\to\infty}n^{\tau D/2-1}L^{-\tau/2}(n)(\log\log n)^{-3\tau/2}\sup_{0\le t\le 1}\sup_{0\le s\le 1}|[nt]Q_n(s,t)|$$
$$=2^{(3\tau+5)/2}\kappa_{3\tau}(\tau!(2-\tau D)(1-\tau D))^{-3/2}\qquad a.s.,$$

*and, as* $n \to \infty$,

$$n^{\tau D/2-1}L^{-\tau/2}(n)\sup_{0\le t\le 1}\sup_{0\le s\le 1}|[nt]Q_n(s,t)|$$
$$\xrightarrow{\mathcal{D}} 2^{5/2}\kappa_{3\tau}((2-\tau D)(1-\tau D))^{-3/2}\sup_{0\le t\le 1}|Y_\tau^3(t)|.$$

*Moreover,*

$$\limsup_{n\to\infty} n^{\tau D/2-1}L^{-\tau/2}(n)(\log\log n)^{-3\tau/2}\|[nt]Q_n\|_p$$
$$= 2^{(3\tau+5)/2}\|J_\tau^2 J_\tau'\|_p(\tau!(2-\tau D)(1-\tau D))^{-3/2} \qquad a.s.,$$

*and, as* $n \to \infty$,

$$n^{\tau D/2-1}L^{-\tau/2}(n)\|[nt]Q_n\|_p \xrightarrow{\mathcal{D}} 2^{5/2}\|J_\tau^2 J_\tau'\|_p((2-\tau D)(1-\tau D))^{-3/2}\|Y_\tau^3\|_p,$$

*where, in both cases,* $Y_\tau$ *is as in Theorem* A.

Reading Theorems 3.2 and 3.3 with $t=1$, they should be compared to Theorem 2.1 and Corollaries 2.1 and 2.2 of [4] in the i.i.d. case.

**4. Sequential general Bahadur–Kiefer processes, strong approximations.** In this section we study the sequential general Bahadur–Kiefer process $R_n(y,t)$ in terms of the sequential uniform Bahadur–Kiefer process $R_n^*(y,t)$.

The following Csáki-type law of the iterated logarithm (cf. [3]) for the sequential uniform quantile process plays a crucial role in comparing the two processes $\rho_n(y,t)$ and $u_n(y,t)$.

PROPOSITION 4.1. *Assume that the assumptions of Proposition* 2.2 *hold. Then, as* $n \to \infty$ *we have*

$$\sup_{\delta_n\le y\le 1-\delta_n}|u_n(y,1)|^2/y(1-y) = O((\log\log n)^\tau) \qquad a.s.,$$

*where* $\delta_n = (n^{-D}L(n)\log\log n)^\tau$.

PROOF. Note that

$$\sup_{\delta_n\le y\le 1/2}|u_n^2(y,1) - |d_n^{-1}V(y,n)|^2|/y$$
$$\le \sup_{0\le y\le 1}|u_n(y,1) - d_n^{-1}V(y,n)|^2 \cdot \delta_n^{-1}$$
$$+ 2\sup_{\delta_n\le y\le 1/2}|d_n^{-1}V(y,n)|/y^{1/2}\sup_{0\le y\le 1}|u_n(y,1) - d_n^{-1}V(y,n)| \cdot \delta_n^{-1/2}.$$

Assumption A and simple calculations yield

(4.1)
$$\sup_{\delta_n \leq y \leq 1/2} |J_\tau(y)|/y^{1/2} = O(1) \quad \text{and}$$
$$\sup_{1/2 \leq y \leq 1-\delta_n} |J_\tau(y)|/(1-y)^{1/2} = O(1)$$

for large enough $n$. Consequently, (2.4), (2.7) and (4.1) imply that, as $n \to \infty$,

$$\sup_{\delta_n \leq y \leq 1/2} |u_n^2(y,1) - |d_n^{-1}V(y,n)|^2|/y = O((\log\log n)^\tau) \quad \text{a.s.}$$

Similarly, as $n \to \infty$ we get

$$\sup_{1/2 \leq y \leq 1-\delta_n} |u_n^2(y,1) - |d_n^{-1}V(y,n)|^2|/(1-y) = O((\log\log n)^\tau) \quad \text{a.s.}$$

This, in turn, results in

(4.2) $$\sup_{\delta_n \leq y \leq 1-\delta_n} \frac{|u_n^2(y,1) - |d_n^{-1}V(y,n)|^2|}{y(1-y)} = O((\log\log n)^\tau) \quad \text{a.s.}$$

On the other hand, by (2.8) and (4.1) we know that, as $n \to \infty$,

$$\sup_{\delta_n \leq y \leq 1-\delta_n} |d_n^{-1}V(y,n)|^2/(y(1-y)) = O((\log\log n)^\tau) \quad \text{a.s.}$$

Thus, via (4.2), as $n \to \infty$ we arrive at

(4.3) $$\sup_{\delta_n \leq y \leq 1-\delta_n} |u_n(y,1)|^2/(y(1-y)) = O((\log\log n)^\tau) \quad \text{a.s.} \qquad \square$$

In light of Proposition 4.1, and Lemma 1 of [7] (cf. Lemma 4.5.2 in [8]), it is natural to introduce the following conditions:

(i) $F$ is twice differentiable on $(a,b)$, where

$$a = \sup\{x : F(x) = 0\}, \qquad b = \inf\{x : F(x) = 1\}, \qquad -\infty \leq a < b \leq +\infty;$$

(ii) $F'(x) = f(x) > 0$ on $(a,b)$;
(iii) for some $0 < \gamma < 1 + (\tau D)/(2 - 2\tau D)$, we have

$$\sup_{a < x < b} F(x)(1-F(x))\frac{|f'(x)|}{f^2(x)} = \sup_{0 < y < 1} y(1-y)\frac{|f'(Q(y))|}{f^2(Q(y))} \leq \gamma;$$

(iv) $A := \overline{\lim}_{x \downarrow a} f(x) < \infty$, $B := \overline{\lim}_{x \uparrow b} f(x) < \infty$;
(v) $\min(A, B) > 0$, or
(v′) if $A = 0$ (resp. $B = 0$), then $f$ is nondecreasing (resp. nonincreasing) on an interval to the right of $a$ (resp. to the left of $b$).

The following proposition concludes a strong approximation of the general quantile process $\rho_n(\cdot,\cdot)$ by $V(\cdot,\cdot)$ of (2.3). Thus, it parallels Proposition 2.2 concerning $u_n(\cdot,\cdot)$, and it is achieved by studying the sup-norm distance between $\rho_n(y,t)$ and $u_n(y,t)$.

PROPOSITION 4.2. *Assume the conditions* (i)–(iii) *on $F$ and the assumptions of Proposition* 2.2. *Then, as $n \to \infty$ we have*

$$\sup_{0 \leq t \leq 1} \sup_{\delta_n \leq y \leq 1 - \delta_n} |\rho_n(y,t) - d_n^{-1} V(y, nt)| \tag{4.4}$$
$$= O(n^{-\tau D/2} L^{\tau/2}(n) (\log \log n)^\tau) \quad a.s.,$$

*where $\delta_n = (n^{-D} L(n) \log \log n)^\tau$. If, in addition to* (i)–(iii), *we also assume* (iv) *and* (v) [*or* (v')], *then*

$$\sup_{0 \leq t \leq 1} \sup_{0 \leq y \leq 1} |\rho_n(y,t) - d_n^{-1} V(y, nt)|$$

$$(4.5) \quad = \begin{cases} O(n^{-\tau D/2} L^{\tau/2}(n) (\log \log n)^{\tau+1}), \\ \quad 0 < \gamma \leq 1, \\ O(n^{(1-\tau D)\gamma + \tau D/2 - 1} L^{\tau \gamma - \tau/2}(n) (\log n)^{(1+C)(\gamma-1)}), \quad a.s., \\ \quad 1 < \gamma < 1 + \dfrac{\tau D}{2(1 - \tau D)}, \end{cases}$$

*where $C > 0$ is arbitrary.*

PROOF. Observe that a two-term Taylor expansion gives

$$\rho_n(y,1) = d_n^{-1} n f(Q(y))(Q(y) - \widehat{Q}_n(y))$$
$$(4.6) \qquad = d_n^{-1} n f(Q(y))(Q(y) - Q(\widehat{U}_n(y)))$$
$$= u_n(y,1) - \frac{d_n}{2n} u_n^2(y) \frac{f'(Q(\theta_{3n}(y)))}{f^3(Q(\theta_{3n}(y)))} f(Q(y)),$$

where $|y - \theta_{3n}(y)| \leq |y - \widehat{U}_n(y)|$.

By (4.1), arguing as in the proof of Theorem 4.5.6 of [8], we arrive at

$$\sup_{0 < \theta_{3n}(y) < 1} \theta_{3n}(y)(1 - \theta_{3n}(y)) \frac{|f'(Q(\theta_{3n}(y)))|}{f^2(Q(\theta_{3n}(y)))} \leq \gamma$$

and

$$\sup_{\delta_n \leq y \leq 1 - \delta_n} \frac{f(Q(y))}{f(Q(\theta_{3n}(y)))} \leq \sup_{\delta_n \leq y \leq 1 - \delta_n} \left[ \frac{\theta_{3n}(y)(1-y)}{y(1 - \theta_{3n}(y))} + \frac{y(1 - \theta_{3n}(y))}{\theta_{3n}(y)(1-y)} \right]^\gamma < \infty.$$

These, together with Proposition 4.1 and (4.6), yield

$$\sup_{\delta_n \leq y \leq 1 - \delta_n} |\rho_n(y,1) - u_n(y,1)| = O(n^{-\tau D/2} L^{\tau/2}(n) (\log \log n)^\tau) \quad \text{a.s.}$$

Hence, arguing as at the end of the proof of Proposition 2.1, we conclude

$$\sup_{0\leq t\leq 1}\sup_{\delta_n\leq y\leq 1-\delta_n}|\rho_n(y,t)-u_n(y,t)|$$
(4.7)
$$= O(n^{-\tau D/2}L^{\tau/2}(n)(\log\log n)^{\tau}) \quad \text{a.s.}$$

Now (2.4) and (4.8) together imply (4.4).

Next, assuming now (iv) and (v), consider the one-term Taylor expansion as in (1.6),

$$\rho_n(y,t) = u_n(y,t)\frac{f(Q(y))}{f(Q(\theta_n(y,t)))}.$$

It follows from Assumption A in combination with (2.4) and (2.7) that

(4.8) $\quad \sup_{0\leq t\leq 1}\sup_{0\leq y\leq \delta_n}|u_n(y,t)| = O(n^{-\tau D}L^{\tau}(n)(\log\log n)^{3\tau/2}) \quad$ a.s.

and

$$\sup_{0\leq t\leq 1}\sup_{1-\delta_n\leq y\leq 1}|u_n(y,t)| = O(n^{-\tau D}L^{\tau}(n)(\log\log n)^{3\tau/2}) \quad \text{a.s.}$$

Hence, we have

$$\sup_{0\leq t\leq 1}\sup_{0\leq y\leq \delta_n}|\rho_n(y,t)-u_n(y,t)|$$
(4.9)
$$= O(n^{-\tau D/2}L^{\tau/2}(n)(\log\log n)^{\tau}) \quad \text{a.s.}$$

and

$$\sup_{0\leq t\leq 1}\sup_{1-\delta_n\leq y\leq 1}|\rho_n(y,t)-u_n(y,t)|$$
(4.10)
$$= O(n^{-\tau D/2}L^{\tau/2}(n)(\log\log n)^{\tau}) \quad \text{a.s.}$$

Using (2.4), (4.8), (4.9) and (4.10), we get

$$\sup_{0\leq t\leq 1}\sup_{0\leq y\leq 1}|\rho_n(y,t)-d_n^{-1}V(y,nt)|$$
(4.11)
$$= O(n^{-\tau D/2}L^{\tau/2}(n)(\log\log n)^{\tau}) \quad \text{a.s.}$$

Finally, we assume (iv) and (v'). In order to prove (4.5), it again suffices to show that $\sup_{0\leq t\leq 1}\sup_{0\leq y\leq \delta_n}|\rho_n(y,t) - u_n(y,t)|$ and $\sup_{0\leq t\leq 1}\sup_{1-\delta_n\leq y\leq 1}|\rho_n(y,t) - u_n(y,t)|$ converge to zero a.s. under assumptions (iv) and (v'). We demonstrate this only for the first one of these, since, for the second one, a similar argument holds.

Along similar lines to the proof of Theorem 4.5.6 in [8], we conclude

$$|\rho_n(y,1)| \leq |u_n(y,1)| \quad \text{if } \widehat{U}_n(y) \geq y,$$

and if $\widehat{U_n}(y) < y$, then

$$|\rho_n(y,1)| \leq \begin{cases} O(n^{\tau D/2} L^{-\tau/2}(n)\delta_n), \\ \quad 0 < \gamma < 1, \\ O(n^{\tau D/2} L^{-\tau/2}(n)\delta_n \log \log n), \\ \quad \gamma = 1 \\ O(n^{\tau D/2} L^{-\tau/2}(n)\delta_n^\gamma n^{\gamma-1}(\log n)^{(1+C)(\gamma-1)}), \\ \quad 1 < \gamma < 1 + \dfrac{\tau D}{2(1-\tau D)}, \end{cases} \quad \text{a.s.,}$$

where $C > 0$ is arbitrary. Note that $-\tau D/2 < (1-\tau D)\gamma + \tau D/2 - 1 < 0$ if $1 < \gamma < 1 + (\tau D)/(2 - 2\tau D)$. Hence, with the help of (4.8), we obtain

$$\sup_{0 \leq y \leq \delta_n} |\rho_n(y,1) - u_n(y,1)|$$

$$= \begin{cases} O(n^{-\tau D/2} L^{\tau/2}(n)(\log \log n)^{\tau+1}), \\ \quad 0 < \gamma \leq 1, \\ O(n^{(1-\tau D)\gamma + \tau D/2 - 1} L^{\tau\gamma - \tau/2}(n)(\log n)^{(1+C)(\gamma-1)}) \\ \quad 1 < \gamma < 1 + \dfrac{\tau D}{2(1-\tau D)}. \end{cases} \quad \text{a.s.,}$$

This, combined with Proposition 2.2 and (4.11), completes the proof of Proposition 4.2. □

REMARK 4.1. Note that

(4.12)
$$\begin{aligned} \{R_n(y,t) - R_n^*(y,t), 0 \leq y, t \leq 1, n = 1, 2, \dots\} \\ = \{-d_n(\rho_n(y,t) - u_n(y,t)), 0 \leq y, t \leq 1, n = 1, 2, \dots\}. \end{aligned}$$

The relationship (4.12) clearly indicates that the results we have summarized and proved in Theorems 2.3–2.5 for $R_n^*(y,t)$ can be immediately restated for the sequential general Bahadur–Kiefer process $R_n(y,t)$ via the strong invariance principle of Proposition 4.2. So we spell out and summarize these results for $R_n(y,t)$ without proof.

THEOREM 4.1. *Assume the conditions* (i)–(iii) *on $F$ and the assumptions of Proposition* 2.2. *Then, as $n \to \infty$ we have*

$$n^{\tau D - 2} L^{-\tau}(n)[nt] R_n(y,t) I\{\delta_n \leq y \leq 1 - \delta_n\}$$

$$\xrightarrow{\mathcal{D}} \frac{2}{(2-\tau D)(1-\tau D)} J_\tau(y) J_\tau'(y) Y_\tau^2(t)$$

*in the space $D[0,1]^2$ equipped with the sup-norm, where $\delta_n = (n^{-D} L(n) \times \log \log n)^\tau$. Moreover,*

$$\limsup_{n \to \infty} n^{\tau D - 2} (L(n) \log \log n)^{-\tau} \sup_{0 \leq t \leq 1} \sup_{\delta_n \leq y \leq 1 - \delta_n} |[nt] R_n(y,t)|$$

$$= \frac{2^{\tau+1}\kappa_{2\tau}}{\tau!(2-\tau D)(1-\tau D)} \qquad a.s.,$$

and, as $n \to \infty$,

$$n^{\tau D-2}L^{-\tau}(n)\sup_{0\le t\le 1}\sup_{\delta_n\le y\le 1-\delta_n}|[nt]R_n(y,t)| \xrightarrow{\mathcal{D}} \frac{2\kappa_{2\tau}}{(2-\tau D)(1-\tau D)}\sup_{0\le t\le 1}Y_\tau^2(t).$$

THEOREM 4.2. *In addition to the conditions in Theorem* 4.1, *we assume* (iv) *and* (v) [*or* (v$'$)]. *Then, as* $n \to \infty$ *we have*

$$n^{\tau D-2}L^{-\tau}(n)[nt]R_n(y,t) \xrightarrow{\mathcal{D}} \frac{2}{(2-\tau D)(1-\tau D)}J_\tau(y)J_\tau'(y)Y_\tau^2(t)$$

*in the space* $D[0,1]^2$, *equipped with the sup-norm, as well as*

$$\limsup_{n\to\infty} n^{\tau D-2}(L(n)\log\log n)^{-\tau}\sup_{0\le t\le 1}\sup_{0\le y\le 1}|[nt]R_n(y,t)|$$

$$= \frac{2^{\tau+1}\kappa_{2\tau}}{\tau!(2-\tau D)(1-\tau D)} \qquad a.s.$$

*and*

$$n^{\tau D-2}L^{-\tau}(n)\sup_{0\le t\le 1}\sup_{0\le y\le 1}|[nt]R_n(y,t)|$$

$$\xrightarrow{\mathcal{D}} \frac{2\kappa_{2\tau}}{(2-\tau D)(1-\tau D)}\sup_{0\le t\le 1}Y_\tau^2(t), \qquad n \to \infty.$$

*Moreover,*

$$\limsup_{n\to\infty} n^{\tau D-2}(L(n)\log\log n)^{-\tau}\|[nt]R_n\|_p = \frac{2^{\tau+1}\|J_\tau J_\tau'\|_p}{\tau!(2-\tau D)(1-\tau D)} \qquad a.s.$$

*and, as* $n \to \infty$,

$$n^{\tau D-2}L^{-\tau}(n)\|[nt]R_n\|_p \xrightarrow{\mathcal{D}} \frac{2\|J_\tau J_\tau'\|_p}{(2-\tau D)(1-\tau D)}\|Y_\tau^2\|_p.$$

**Acknowledgments.** We sincerely wish to thank two referees and an Associate Editor for their queries and many insightful remarks and suggestions which have led to improving the presentation of our results, and to making corrections of some oversights in our proofs. A question posed by the Associate Editor resulted in removing an unnatural condition from the statement of Proposition 2.1.

M. CSÖRGŐ  
B. SZYSZKOWICZ  
SCHOOL OF MATHEMATICS  
AND STATISTICS  
CARLETON UNIVERSITY  
1125 COLONEL BY DRIVE  
OTTAWA, ONTARIO  
CANADA K1S 5B6  
E-MAIL: mcsorgo@math.carleton.ca  
bszyszko@math.carleton.ca  

L. WANG  
DEPARTMENT OF MATHEMATICS  
NANJING UNIVERSITY  
NANJING 210093  
P.R. CHINA